\documentclass{amsart}
\usepackage{graphicx}
\theoremstyle{definition}

\theoremstyle{remark}

\numberwithin{equation}{section}

\begin{document}

\title{Uniform approximation of some Dirichlet series by partial products of Euler type}
\author{I.Sh.Jabbarov}
\address{AZ2000 Haydar Aliyev avenue, 187, Ganja State University, Azerbaijan}
\email{jabbarovish@rambler.ru}


\dedicatory{Dedicated to the memory of professor Voronin S. M.}
\begin{abstract}
In the work it is gotten a uniform approximation of Dirichlet series
defined by Euler product by patrial products of Euler type in discs
placed in the right half of the critical strip. As a consequence the
analog of the Riemann Hypothesis is proven.
\end{abstract}
\maketitle
\section{Introduction.}

Appearance of Dirichlet series and understanding of their fundamental role in Analytical Number Theory is connected with L.Euler's name. In 1748 Euler
([1]) entered the zeta-function $\zeta (s)=\sum n^{-s}  $considering it as a function of real variable s, and proved an important identity:
\[\zeta (s)=\prod _{p} \left(1-p^{-s} \right)^{-1} ,s>1;\]
here the product is taken over all prime numbers. This product is called to be Euler product.

	In 1837 using and developing Euler's ideas, L.Dirichlet gave generalization of the theorem of Euclid for arithmetic progressions considering L
-functions. Investigations of Dirichlet showed an importance of studying of Dirichlet series defined by more general Euler products.

	After coming into the world of essential Riemann's work [2] it stood clear that deeper results of the theory of prime numbers are connected with
analytical properties of Dirichlet series of complex variable.

    Many of analytical properties of Dirichlet series were studied by using of various their finite approximations. For example, some questions of the
theory of Dirichlet series connected with mean values, an order, or a density distribution of zeroes investigated in the works [3-8] by using of
approximations by partial sums. G. Bohr and E. Landau were the first who applied partial products of Euler type to investigate the zero distribution of
the zeta function ([6-8]). In the works [9-16] S. M. Voronin developing the method of G. Bohr and E. Landau used special type of partial products in
the questions of distribution of zeroes and non-zero values of L-functions in the critical strip. In compliance with the Universality Theory of S. M.
Voronin every analytical function non vanishing in and on the disc $\left|s\right|=r<1/4$, can be approximated by finite products of the form
$\prod _{p}\left(1-e^{-2\pi i\vartheta _{p} } p^{-s} \right) ^{-1} $ where \textit{p} takes on values from some finite set of prime numbers.

In the present work we show that the Dirichlet series with the Euler product having analytical continuation to the critical strip without
singularities in some natural conditions can be approximated by partial products of Euler type in the discs of the critical strip where the primes,
over which the products are taken, are distributed by a suitable way (see formulation of the theorem below) (see [29,30]). The family of such series
includes many of widely used Dirichlet series as the zeta function, Dirichlet L - functions, or L-functions of some algebraic extensions with the
commutative Galois groups and etc.

 Let we are given with a following infinite product taken over all prime numbers \textit{p}:
\begin{equation} \label{1}
F(s)=\prod _{p}f_{p} (p^{-s} ) ,
\end{equation}
where $f_{p} (z)$ is a rational function of a variable \textit{z} having not poles in the disc $\left|z\right|<1$,
\[f_{p} (z)=1+\sum _{m=1}^{\infty }a_{p}^{m}  z^{m} ,\]
and for any positive small $\varepsilon $ the inequality $\left|a_{p}^{m} \right|\le c(\varepsilon )p^{m\varepsilon } ;c(\varepsilon )\ge 1$ is
satisfied uniformly by \textit{p}.

\textbf{\textit{Theorem. }}\textit{Let the function }$F(s)$ \textit{ have not singularities in the half plane $\sigma >1/2$ with exception of
finite number of poles on the line$\sigma =1$, and every factor of the product (1) have not zeroes in the half plane }$\sigma >1/2$. \textit{Suppose that
for any small positive number $\lambda $ there exist constants $c_{0} =c_{0} \left(\lambda \right)>0$  and $h_{0} =h_{0} \left(\lambda \right)>0$,
satisfying, for any $h>h_{0} $, the following inequality:}
\[\sum _{h<p\le h(1+\log ^{-10} h)}\left|a_{p}^{1} \right| p^{-(1-\lambda )} \ge c_{0} (\lambda )h^{\lambda /4} .       \]

\textit{If in the disc }$\left|s-\sigma _{0} \right|\le r<r_{0} =\min (1-\sigma _{0} ,\sigma _{0} -1/2)$ $F\left(s+it_{0} \right)$ \textit{ has not zeros for some real $t_{0} $, and $1/2<\sigma _{0} <1$ then there exists a  sequence $(\theta _{n} )$, $\theta _{n} \in \Omega =[0,1]\times [0,1]\times
\cdots $ and a sequence $(m_{n} )$ of integers that }
\[\mathop{\lim }\limits_{n\to \infty } F_{n} (s+it,\theta _{n} )=F(s+it+it_{0} ),\]
\textit{for every real t uniformly by s in this disc; here }$\theta _{n} =(\theta _{p}^{n} )$\textit{,  and}
\begin{equation} \label{2}
F_{n} (s+it,\theta _{n} )=\prod _{p\le m_{n} }f_{p} (e^{-2\pi i(\theta _{p}^{n} +\gamma _{p} )} p^{-s-it} ) ;\; \gamma _{p} =\frac{t_{0} \log p}{2\pi
}.
\end{equation}
If we put now $a_{p}^{m} =1$ for all natural $m$ and prime $p$ we get the zeta function. From the theorem of Valle-Poussin C. J. [see, 4, p. 59] it
follows that the all conditions of the theorem formulated above are satisfied for the Riemann zeta function with $\sigma _{0} =3/4$ and $t_{0} =0$.
Following by [30] we prove the analog of the Riemann Hypothesis for the function $F(s)$.

\textbf{\textit{Corollary. }}\textit{The analog of the Riemann Hypothesis is true: F(s)} $\ne $\textit{0, when }$1/2<\sigma <1$\textit{. }

\section{Additional statements.}

\textbf{Lemma 1.} \textit{Let a series of analytical functions }

\[\sum _{n=1}^{\infty } f_{n} (s)\]
\textit{be given in one-connected domain }$G$\textit{ of a complex }$s$\textit{-plane and converges absolutely almost everywhere in }$G$\textit{ in the
Lebesgue sense and the function }
\[\Phi (\sigma ,t)=\sum _{n=1}^{\infty } |f_{n} (s)|\]
\textit{is a summable function in }$G$. \textit{ Then, the given series converges uniformly in any compact subdomain of }$G$\textit{; particularly,
the sum of this series is an analytical function in }$G$.
\textit{Proof}. It is enough to show that the theorem is true for any rectangular area in $G$. Let $C$ be a rectangle in $G$ and $C'$ another rectangle
inside of $C$ and their sides are parallel to the co-ordinate axes. We can assume that on a contour of these rectangles the series converges almost
everywhere, according to the theorem of Fubini (see [17, p. 208]). Let $\Phi _{0} (s)=\Phi _{0} (\sigma ,t)$ be the sum of given series at points
of convergence. Under the theorem of Lebesgue on a bounded convergence (see [20, p. 293]), we have:
\[(2\pi i)^{-1} \int _{C} \frac{\Phi _{0} (s)}{s-\xi } ds=\sum _{n=1}^{\infty } (2\pi i)^{-1} \int _{C} \frac{f_{n} (s)}{s-\xi } ds,\]
where the integrals are taken in Lebesgue sense. As on the right part of the last equality the integrals are existing in the Riemann sense also, then by applying Couchy's formula and denoting the left side of the last equality by $\Phi _{1} (\xi )$  we get for any point $\xi $ on or in the contour
of $C'$:
\[\Phi _{1} (\xi )=(2\pi i)^{-1} \int _{C} \frac{\Phi _{0} (s)}{s-\xi } ds=\sum _{n=1}^{\infty } f_{n} (s)\]
(certainly, $\Phi _{1} (\xi )=\Phi _{0} (\xi )$ almost everywhere). Further, the series can be estimated in $C'$ by a following way
\[|f_{n} (\xi )|\le (2\pi )^{-1} \int _{C} \frac{|f_{n} (s)|}{|s-\xi |} |ds|\le (2\pi \delta )^{-1} \int _{C} |f_{n} (s)||ds|,\]
when $\delta $ designates the minimum of the distances between the sides of $C$and$C'$. The series
\[\sum _{n=1}^{\infty } \int _{C} |f_{n} (s)||ds|\]
converges, in the consent with the theorem of Lebesgue on monotone convergence (see [20, p. 290]). Hence, the series $\sum _{n=1}^{\infty } f_{n} (\xi
)$ converges uniformly in the inside of $C'$. The lemma 1 is proved.

Let's following by authors of [4] enter the notion of Hardy space (see also [30]).

\textbf{Definition 1. }\textit{The set of functions} $f(s)$ \textit{ defined for }$\left|s\right|<R$\ \textit{ and analytical in this area, for which}
\[\left\| f\right\| ^{2} =\mathop{\lim }\limits_{r\to R} \int\int_{\left|s\right|<r}\left|f(s)\right| ^{2} d\sigma dt<\infty ;s=\sigma +it,\]
\textit{is called the Hardy space and is designated as} $H_{2}^{(R)}, R>0$.

Obviously, Hardy space is a real linear space in which is possible to enter a scalar product of functions by means of  the equality
\begin{equation} \label{3}
\left(f(s),g(s)\right)=Re\int\int_{\left|s\right|\le R}f(s)\overline{g(s)} d\sigma dt.
\end{equation}

Using the entered scalar product, we will prove that Hardy space is a real Hilbert space.

  \textbf{Lemma 2. }\textit{The Hardy space} $H_{2}^{(R)} $\textit{ together with the entered scalar product (3) is a real Hilbert space.}

\textit{Proof}. It is enough to prove that any fundamental sequence $\left(f_{n} (s)\right)_{n\ge 1} $ converges to some analytical function
$f\left(s\right)\in H_{2}^{\left(R\right)} $. As the sequence is fundamental there exist such a sequence of natural numbers $\left(n_{k} \right)_{k\ge 0} $ that for any natural $k$ we have:
\[\left\| f_{n_{k} } -f_{n_{k-1} } \right\| \le 2^{-k} .\]
Let's consider a series of analytical functions
\[f_{n_{0} } +\sum _{k=1}^{\infty }\left(f_{n_{k} } -f_{n_{k-1} } \right) .\]
We will prove first that it converges uniformly in any closed disc lying in the open disc $\left|s\right|<R$. According to the definition of the norm
we have:
\[\left\| f(s)\right\| ^{2} =\int_{\left|s\right|<R}\left|f(s)\right|^{2} d\sigma dt ,\]
possible, in improper meaning of definition of the norm. Then, designating
\[g(s)=\sum _{k=1}^{\infty }\left|f_{n_{k} } (s)-f_{n_{k-1} } (s)\right| \]
we receive:
\[\int_{\left|s\right|<R}g(s)d\sigma dt\le \sum _{k=1}^{\infty }\left(\pi R^{2} \iint_{\left|s\right|<R}\left|f_{n_{k} } -f_{n_{k-1} } \right|^{2}
d\sigma dt \right)^{1/2} \]
\[  \le \sqrt{\pi } R\sum _{k=1}^{\infty }2^{-k} <+\infty  .\]
Hence, the function $g(s)$ is a summable function in variables $\sigma ,t$ then the lemma 1 is applicable. Applying the lemma 1 we receive  that the
series $f_{n_{0} } +\sum _{k=1}^{\infty }$ $\left(f_{n_{k} }-f_{n_{k-1} } \right) $ converges uniformly in any disc $\left|s\right|\le r<R$. Then the
subsequence $\left(f_{n_{k} } (s)\right)_{k\ge 1} $ converges to some analytical function $\varphi (s)$. As the sequence is fundamental then for any
$\varepsilon >0$ it can be found $n_{0} $ such that for any natural $m>n_{0} $ the inequality
\[\iint_{\left|s\right|<R}\left|\varphi (s)-f_{m} (s)\right|^{2} d\sigma dt<\varepsilon  \]
holds. Let $r<R$ be any real number. Then, using an inequality of [19, p. 345] one can receive
\[r^{2} \left|\varphi (s)-f_{m} (s)\right|^{2} \le \pi ^{-1} \iint_{\left|s\right|<R}\left|\varphi (s)-f_{m} (s)\right|^{2} d\sigma dt<\varepsilon /\pi
 \]
for any $s,\left|s\right|\le r$. As $\varepsilon $ is arbitrarily then from here it follows the convergence of the sequence $\left(f_{m}
(s)\right)_{m\ge 1} $ to $\varphi (s)$. As
\[\iint_{\left|s\right|<R}\left|\varphi (s)\right|^{2} d\sigma dt\le  \iint_{\left|s\right|<R}\left|f_{n_{0} } (s)\right|^{2} d\sigma dt+
\iint_{\left|s\right|<R}\left|g(s)\right|^{2} d\sigma dt<+\infty  \]
then $\varphi \left(s\right)\in H_{2}^{\left(R\right)}, $ and therefore, the considered space is complete. The lemma 2 is proved.

The following is a variant of S. M. Voronin's lemma proved by him in [4] for the zeta function, and it is one of the basic arguments of our work.

\textbf{\textit{Lemma 3. }}\textit{Let }$0<r<1/4$, \textit{and }$g(s)$ \textit{is an analytical function in the cdisc}$\left|s\right|\le r$\textit{,
being continuous and does not vanishing on the circle $\left|s\right|= r$. Then, for any }$\varepsilon >0$\textit{ and }$y>2$\textit{ there exist
a finite set of prime numbers M containing all of primes }$p$, $p\le y$ \textit{ for which the following inequality is fair:}
\[\mathop{\max }\limits_{\left|s\right|\le r} \left|g(s)-F_{M} \left(s+\frac{3}{4} ;\theta \right)\right|\le \varepsilon \]
\textit{for some} $\theta =(\theta _{p} )_{p\in M} $\textit{ with} $\theta _{p} =\theta _{p}^{0} $ \textit{being numbers set beforehand from the interval } $\left[0,1\right]$\textit{ for }$p\le y$\textit{; the function }$F_{M} \left(s+\frac{3}{4} \; \theta \right)$\textit{ is given by the
equality}
\[F_{M} (s+3/4;\theta )=\prod _{p\in M}f_{p} \left(e^{-2\pi i\theta _{p} } p^{-s-3/4} \right) .\]

\textit{   Proof. }The proof of the lemma 3 we will spend by the method of the work [4] of  Voronin S. M. As $g(s)$ is an analytical function in the
circle $\left|s\right|\le r$ then we will consider an auxiliary function $g(s/\gamma ^{2} )$ $(\gamma >1,\gamma ^{2} r<1/4)$ which for any $\varepsilon
>0$ at some $\gamma $ satisfies the inequality $\mathop{\max }\limits_{\left|s\right|\le r} \left|g(s)-g(s/\gamma ^{2} )\right|<\varepsilon $.
Therefore, it is enough to prove the statement of the lemma 3 for the function $g(s/\gamma ^{2} )$ in the disc $|s|\le r$. The advantage is consisted in that that the function $g(s/\gamma ^{2} )$ belongs to the space $H_{2}^{(\gamma r)} $ (a disc has a radius greater than $r$ which is important for our subsequent reasoning). Not breaking, therefore, a generality we believe that the function $g(s)$ is an analytical in the disc $|s|\le r\gamma ^{2}$ and we will consider the space $H_{2}^{(\gamma r)} $.

The function $\log g(s)$ in the conditions of the theorem has no singularities in the disc $|s|\le r\gamma $. Therefore, it is enough to prove an existence of a such element $\theta $, satisfying the conditions of the lemma 3, that
\[\mathop{\max }\nolimits_{|s|\le r} |\log g(s)-\log F_{M} (s+3/4;\theta )|\le \varepsilon .\]
A series (2) of the work [4, p.241] we define as
\[u_{k} (s)=\log f_{p} \left(e^{-2\pi i(\varphi _{k} +\vartheta _{k} )} p_{k} ^{-s-3/4} \right),   \]
supposing $2\pi \varphi _{k} $ to be an argument of the coefficient $a_{p_{k} }^{1} $. We have:
\[u_{k} (s)=\log (1+a_{p_{k} }^{1} e^{-2\pi i(\varphi _{k} +\vartheta _{k} )} p_{k} ^{-s-3/4} )+\log (1+a_{p_{k} }^{1} e^{-2\pi i(\varphi _{k}
+\vartheta _{k} )} p_{k} ^{-s-3/4} )^{-1}\times\]
\[ f_{p} (e^{-2\pi i(\varphi _{k} +\vartheta _{k} )} p_{k} ^{-s-3/4} ), \]
and using decomposition of the logarithmic function into power series, we get
\begin{equation} \label{4}
u_{k} (s)=a_{p_{k} }^{1} e^{-2\pi i(\varphi _{k} +\vartheta _{k} )} p_{k}^{-s-3/4} + \nu_{k} (s),
\end{equation}
for every $k$ being large enough. So, for any $\varepsilon >0$
\[\nu_{k} (s)=O(p_{k}^{2\varepsilon +2r-3/2} )+\log \left(1+\sum _{m=2}^{\infty }b_{m} (e^{-2\pi i(\varphi _{k} +\vartheta _{k} )} p_{k} ^{-s-3/4}
)^{m}\right),\]
and the factors $b_{n} $ are defined by a following equality:
\[b_{m} =a_{p}^{m} -a_{p}^{m-1} a_{p}^{1} +a_{p}^{m-2} \left(a_{p}^{1} \right)^{2} -\cdots +(-1)^{m-2} a_{p}^{2} \left(a_{p}^{1} \right)^{m-2} .\]
We have
\[\left|b_{m} \right|\le (m-1)c^{m-1} (\delta )p^{\varepsilon \, m} .  \]
As $r<1/4$, we can take $\varepsilon >0$ such that the inequality $4\varepsilon +2r-3/2<-1$ was satisfied. Then, definition of $u_{m} (s)$ and (4),
together with the last inequality, show that the series
\begin{equation} \label{5}
\sum _{m=1}^{\infty }\eta _{m} (s) ;\; \eta _{m} (s)=a_{p_{m} }^{1} e^{-2\pi i(\varphi _{m} +\vartheta _{m} )}
\end{equation}
differs from the series $\sum u_{m} (s) $ by an absolutely converging series. Really, since $m\le 2^{m} $ then
\[\sum _{m=2}^{\infty }\left|b_{m} \right|p^{-m(\varepsilon +r-3/4)} \le \sum _{m=2}^{\infty }2^{m-1} c^{m-2} (\varepsilon )p^{m(\varepsilon +r-3/4)}
\le \]
\[\frac{2c(\varepsilon )p^{2\varepsilon +2r-3/2} }{1-2c(\varepsilon )p^{\varepsilon +r-3/4} } \le 4c(\varepsilon )p^{-2\varepsilon -1} ,\]
if  $p$ is so large that $2c(\varepsilon )p^{r-3/4} \le 1/2$. Therefore, the series $\sum \left|\nu _{k} \right| $ converges.

 Now, it is enough for us to show that for any $\varphi (s)\in H_{2}^{\left(\gamma \, r\right)} $ ($0\, <\gamma <1$ is any) there exists some subseries of (5) converging to $\varphi (s)$. In particular, taking $n$ equal to the greatest value of $k$ for which $p_{k} \le y$ we admit
\[\varphi (s)=\log g(s)-\sum _{k>n}(u_{k} (s)-\eta _{k} (s))-\sum _{k\le n}u_{k} (s)  .\]
Considering the last remark we will find some permutation of  $\sum _{k>n}\eta _{k} (s) $ converging to $\varphi (s)$ (clearly, any permutation of the series $\sum _{k>n}(u_{k} (s)-\eta _{k} (s)) $ converges to the same sum uniformly). Then for any $\varepsilon $ there will be found such a set of indexes $k\in M$ that
\[\mathop{\max }\limits_{\left|s\right|\le r} \left|\varphi (s)-\sum _{k\in M,p_{k} >y}\eta _{k} (s) \right|\le \varepsilon /2.\]
Let $q(s)=\sum _{k=n+1}^{\infty }\left(u_{k} (s)-\eta _{k} (s)\right) $. As this series converges absolutely it is possible to select mentioned above set $M$ so that the following relationship was carried out
\[\left|q(s)-\sum _{k\in M,k>n}\left(u_{k} (s)-\eta _{k} (s)\right) \right|\le \varepsilon /2.\]
Then we will receive:
\[\left|\varphi (s)-\sum _{k\in M,p_{k} >y}\eta _{k} (s) \right|=\left|\log g(s)-\sum _{n\in M}u_{n} (s) \right|\le \varepsilon .\]
Thereby, the proof of the lemma 3 will be finished.

At $k>n$ we set $\theta _{k} =\rho (k)/4$ and $\rho (k)$ will be defined below. For $k\le n$ we take $\theta _{k} +\varphi _{k} =0$. Then, for $k>n$, we have:
\[(\eta _{k} (s),\varphi (s))=\left|a_{p_{k} }^{1} \right|Re\int _{\left|s\right|\le R}e^{-2\pi i\rho (k)/4} p_{k} ^{-(s+3/4)} \overline{\varphi (s)}d\sigma dt=\]
\[=Re\left[\left|a_{p_{k} }^{1} \right|e^{-2\pi i\rho (k)/4} \Delta (\log p_{k} )\right] ,\]
by denoting
\[\Delta \left(x\right)=\iint \nolimits _{\left|s\right|\le R}e^{-x(s+3/4)} \overline{\varphi (s)}d\sigma dt. \]

Writing $R=\gamma r$ consider the space $H_{2}^{(R)} $. Then,
\[\left\| \eta _{k} (s)\right\| ^{2} =\iint_{\left|s\right|\le R}\left|e^{-2\pi i\theta _{k} } p_{k}^{-s-3/4} \right|^{2} d\sigma dt\le  \pi R^{2} p_{k}^{2r-3/2} .\]
Hence,
\[\sum _{k=1}^{\infty }\left\| \eta _{k} (s)\right\| ^{2} \le \pi R^{2}  \sum _{k=1}^{\infty }p_{k}^{2r-3/2} <+\infty  ,\]
the first condition of the theorem 1, \S 6 of an appendix of [4] is executed.

Let now $\varphi (s)\in H_{2}^{(R)} $ be arbitrary element of the space with the condition $\left\| \varphi (s)\right\| ^{2} =1$.
Let $\varphi (s)$ have a following expansion into a power series in the disc $\left|s\right|\le R$:
\[\varphi (s)=\sum _{n=0}^{\infty }\alpha _{n} s^{n}  .\]
Then,
\[1=\iint_{\left|s\right|\le R}\left|\sum _{n=0}^{\infty }\alpha _{n} s^{n}  \right| ^{2} d\sigma dt.\]
Exchange the variables under the integral by formulas: $\sigma =r\cos \varphi ,\; t=r\sin \varphi ,\; r\le R$,$0\le \varphi <2\pi $. Then,
\[1=\sum _{n=0}^{\infty }\sum _{n=0}^{\infty }\alpha _{n} \bar{\alpha }_{m}   \int _{0}^{R}r^{n+m+1} \int _{0}^{2\pi }(\cos 2\pi (n-m)\varphi +i\sin 2\pi (n-m)\varphi )  d\varphi .\]
The inner integral is equal to 0 when $m\ne n$, and $2\pi $ otherwise. Hence,
\begin{equation} \label{6}
\pi \sum _{n=0}^{\infty }\left|\alpha _{n} \right|^{2} R^{2n+2}  (n+1)^{-1} =1.
\end{equation}
Let's prove now that there exist a point $\theta $ non dependent on the function $\varphi \left(s\right)$ such that the series $\sum _{k=1}^{\infty }\left(\eta _{k} (s),\varphi (s)\right) $ converges after of some permutation of its members. We have
\[(\eta _{k} (s),\varphi (s))=-Re\int  \int _{|s|\le R} e^{-2\pi i\theta _{k} } p_{k}^{-s-3/4} \overline{\varphi (s)}d\sigma dt=Re[-e^{-2\pi i\theta
_{k} } \Delta (\log p_{k} )].\]
It is possible to represent the function $\Delta (x)$ by a following way:
\[\Delta (x)=e^{-3x/4} \iint_{\left|s\right|\le R}\left(\sum _{n=0}^{\infty }(-sx)^{n} /n! \right) \overline{\left(\sum _{n=0}^{\infty }\alpha _{n} s^{n}  \right)}d\sigma dt=\]
\[=\pi R^{2} e^{-3x/4} \sum _{n=0}^{\infty }(-1)^{n} \bar{\alpha }_{n} x^{n} R^{2n} /(n+1)! =\pi R^{2} e^{-3x/4} \sum _{n=0}^{\infty }\beta _{n}(xR)^{n} /n! \]
by denoting $\beta _{n} =(-1)^{n} R^{n} \bar{\alpha }_{n} /(n+1)$. From (6) one may conclude:
\[\sum _{n=1}^{\infty }\left|\beta _{n} \right| ^{2} \le 1.\]
Hence, $\left|\beta _{n} \right|\le 1$, and, therefore, the function
\begin{equation} \label{7}
H(u)=\sum _{m=0}^{\infty } \frac{\beta _{m} }{m!} u^{m}
\end{equation}
will be an entire function. So,
\[\Delta (x)=\pi R^{2} e^{-3x/4} H(xR).\]
Let's prove that for any $\delta >0$ there will be found tending to the infinity sequence $u_{1} ,u_{2} ,...$, satisfying the inequality
\begin{equation} \label{8}
|H(u_{j} )|>ce^{-(1+2\delta )u_{j} } .
\end{equation}
Let's admit an opposite, i.e. let there exists a positive number $\delta <1$ such that at some $A>0$ being large enough the following inequality
\[\left|H(u)\right|\le Ae^{-(1+2\delta )u} \]
is executed for all $u\ge 0$; in this case we have:
\[\left|e^{(1+\delta )u} H(u)\right|\le Ae^{-\delta \left|u\right|} \; ;u\ge 0.\]
From proved above for $u<0$ one receives:
\[\left|H(u)\right|\le \sum _{n=0}^{\infty }\left|u\right|^{n} /n!=e^{-u}  .\]
Then we have
\[\left|e^{(1+\delta )u} H(u)\right|\le e^{\delta \, u} \le e^{-\delta \left|u\right|} .\]
Consequently, the integral below is existing:
\[\int _{-\infty }^{\infty }\left|e^{(1+\delta )u} H(u)\right|^{2} du .\]
As the function (7) is an entire function of exponential type then the function $e^{(1+\delta )u} H(u)$ will be such one also and belong to the class $E^{\sigma } $ (see [4, p. 408]) with $\sigma <3$. Then under the theorem of Paley -- Wiener (see at the same work) it is existing a finitary
function $h(\xi )\in L_{2} (-3,3)$ such that
\[e^{(1+\delta )u} F(u)=\int _{-3}^{3}h(\xi )e^{iu\xi } d\xi  .\]
Taking converse transformation, we find:
\[h(\xi )=\frac{1}{2\pi } \int _{-\infty }^{\infty }\left(e^{(1+\delta )u} H(u)\right)e^{-iu\xi } du .\]
From the found above estimations it follows that this integral converges absolutely and uniformly in the strip $\left|Im\xi \right|<\delta /2$, and, consequently, represents an analytical function in this strip which contradicts finitaryness of $f(\xi )$. The received contradiction proves an existence of a sequence of points with the condition (8).

Denoting $x_{j} =u_{j} /R$, on the basis of (8) we can assert that
\[\left|\Delta (x_{j} )\right|>ce^{-3x_{j} /4} \left|H(x_{j} R)\right|\ge ce^{-3x_{j} /4} e^{-(1+2\delta )x_{j} R} =ce^{-x_{j} (R+2\delta R+3/4)} .\]
If $\delta >0$ is sufficiently small then $R+2\delta R+3/4<1$, and, hence, there is a $\delta _{0} >0$ such that
\[\left|\Delta (x_{j} )\right|>e^{-(1-\delta _{0} )x_{j} } .   \]
Let's consider the function $\Delta (x)$ on the segment [$x_{j} -1,x_{j} +1$]. Following by [4] we denote $N=\left[x_{j} \right]+1$. From the
estimation for the factors $\beta _{n} $ it follows the inequality:
\[\left|\sum _{n=N^{2} +1}^{\infty }\frac{\beta _{n} }{n!} (xR)^{n}  \right|\le \sum _{n=N^{2} +1}^{\infty }\frac{(xR)^{n} }{n!}  \le
\frac{(xR)^{N^{2}
} }{(N^{2} )!} \sum _{n=0}^{\infty }\frac{(xR)^{n} }{n!} \le \frac{(xR)^{N^{2} } }{(N^{2} )!} e^{N}  .\]
If $n,m\ge 0$ are integers, then $(n+m)!=n!(n+1)\cdots (n+m)\ge n!m!.$ For enough large natural $m$ one has under Stirling's formula:
\[m!=\Gamma (m+1)\ge e^{m\log m-m} =(m/e)^{m} .\]
Hence,
\[\left|\sum _{n=N^{2} +1}^{\infty }\frac{\beta _{n} }{n!} (xR)^{n}  \right|\le \frac{(xR)^{N^{2} } }{(N^{2} )!} e^{N} \le N^{N^{2} }
\left(\frac{N^{2}
}{e} \right)^{N^{2} } e^{N} <<e^{-2x_{j} } \]
at $x\in $ [$x_{j} -1,x_{j} +1$]. Further, $\sum _{n=0}^{N^{2} }\beta _{n} (xR)^{n} /n!<<e^{xR}  $. Analogically,
\[\left|\sum _{n=N^{2} +1}^{\infty }\frac{(-3x/4)^{n} }{n!}  \right|\le \frac{(3x/4)^{N^{2} } }{(N^{2} )!} \sum _{n=0}^{\infty }\frac{(3x/4)^{n} }{n!}
\le \frac{(3x/4)^{N^{2} } }{(N^{2} )!} e^{N} <<e^{-2x_{j} }  ,\]
and $\sum _{n=0}^{N^{2} }(-3x/4)^{n} /n!<<e^{3x/4}  $ for $x\in $ [$x_{j} -1,x_{j} +1$]. Thus,
\[\Delta (x)=\pi R^{2} \sum _{n=0}^{N^{2} }\frac{\left(-3x/4\right)^{n} }{n!}  \sum _{n=0}^{N^{2} }\frac{\beta _{n} }{n!} \left(xR\right)^{n}
+O(e^{-x_{j} } )=\sum _{n=0}^{N^{4} }a_{n} x^{n} +O(e^{-x_{j} } ) \]
According to (8) we receive an inequality
\[\mathop{\max }\limits_{|x-x_{j} |\le 1} |\Delta (x)|>e^{-(1-\delta _{0} )x_{j} } \]
for any $j=1,\, 2,\, ...$. Let $a_{n} =b_{n} +ic_{n} ,b_{n} ,c_{n} \in R$. Then,
\[\Delta (x)=\sum _{n=0}^{N^{4} }b_{n} x^{n} +i \sum _{n=0}^{N^{4} }c_{n} x^{n} +O(e^{x_{j} } ). \]
Therefore, for any $j$ at least one of the following inequalities is executed:
\[\mathop{\max }\limits_{|x-x_{j} |\le 1} \left|\sum _{n=0}^{N^{4} }b_{n} x^{n}  \right|>0.1e^{-(1-\delta _{0} )x_{j} } ,\]
or
\[\mathop{\max }\limits_{|x-x_{j} |\le 1} \left|\sum _{n=0}^{N^{4} }c_{n} x^{n}  \right|>0.1e^{-(1-\delta _{0} )x_{j} } .\]

Let's consider the first possibility. Let $x_{0} $ be the point at which the maximum of modulus is reached. We will designate by $\tau _{j} $ a segment which lies in the interval [$x_{j} -1,x_{j} +1$] containing the point $x_{0} $, and each point $x$ of $\tau _{j} $ satisfies the inequality:
\[\left|g(x)\right|\ge 0.1\left|g(x_{0} )\right|.\]

Let for definiteness $g(x_{0} )<0;\, g(x)=\sum _{n=0}^{N^{4} }b_{n} x^{n}  $. If
\[\tau _{j} \ne  [x_{j} -1,x_{j} +1]\]
(the case of coincidence of intervals is trivial) then there exist a point $x_{1} \in \tau _{j} $ for which
\[\left|g(x_{1} )\right|\le 0.1\left|g(x_{0} )\right|.\]
Now we have:
\[\left|g(x_{0} )-g(x_{1} )\right|\ge 0.5\left|g(x_{0} )\right|.\]
Under the theorem of Lagrange there exist a point $y_{j} \in \tau _{j} $ such that
\[\left|g'(y_{j} )(x_{1} -x_{0} )\right|\ge 0.5\left|g(x_{0} )\right|.\]
Applying the theorem 9, \S 2 of  the appendix of [4] we find:
\[N^{8} \left|g(x_{0} )\right|\left|x_{1} -x_{0} \right|\ge \left|g'(y_{j} )(x_{1} -x_{0} )\right|\ge 0.5\left|g(x_{0} )\right|.\]
So, the interval $\tau _{j} $ has a length not less than $0.5x_{j}^{-8} $. Accepting $h=e^{x_{j} } $, we notice that $[h,h(1+\log ^{-10} h)]\subset$
$[e^{\alpha } ,e^{\alpha +\beta } ]$. From the condition (2) it follows that the set of primes, satisfying the condition $h<p\le h(1+\log ^{-10} h)$ can be distributed among subsets $P_{1} ,P_{2} ,P_{3} ,P_{4} $ for each of which the following inequality is fulfilled
\[\sum _{p\in P_{i} ,\, h<p\le h(1+\log ^{-10} h)}\left|a_{p}^{1} \right| p^{-(1-\lambda )} \ge 0.1c_{0} (\lambda )h^{\lambda /4} ,i=1,2,3,4.\]

To prove the last statement, at first, we divide the set of primes with $h<p\le h(1+\log ^{-10} h)$ into the union of fore subsets $P'_{1} ,P'_{2}
,P'_{3} ,P'_{4} $ arbitrarily. Then, for one of subsets, say for the $P'_{1} $, we will have
\[\sum _{p\in P'_{1} }\left|a_{p}^{1} \right| p^{-(1-\lambda )} \ge 0.25c_{0} (\lambda )h^{\lambda /4} .\]
Now we return the primes corresponding to addends of the last sum back to the union $P'_{2} \bigcup P'_{3} \bigcup P'_{4} $ omitting at the same time addends from this sum consequently, until it is not less than $0.2c_{0} (\lambda )h^{\lambda /4} $. The last returned addend has a bound $\left|a_{p}^{1}
\right|p^{-(1-\lambda )} \le c(\varepsilon )p^{\varepsilon +\lambda -1} $. So, denoting by $P_{1} $ the set remaining after of described above taking off addends for great values of $h_{0} $ we have
\[\sum _{p\in P_{1} }\left|a_{p}^{1} \right| p^{-(1-\lambda )} \ge 0.1c_{0} (\lambda )h^{\lambda /4} .\]
At the same time
\[\sum _{p\in P'}\left|a_{p}^{1} \right| p^{-(1-\lambda )} \ge 0.8c_{0} (\lambda )h^{\lambda /4} ,\]
where the sum over the set $P'$ contains the addends not belonging into the subset$P_{1} $. Continuing the same reasoning we construct suitable subsets.

We put $\rho (k)=i-1$ for every $p_{k} \in P_{i} $. Therefore, if  $\lambda <\delta _{0} /2$ then denoting by $\gamma_{j}$ a corresponding interval of variance of numbers $log p_{k}$, we have:
\[\sum _{\begin{array}{l} {\quad \quad p_{k} \in P_{1} } \\ {\, \log p_{k} \in \gamma _{j} } \end{array}}Re\left[\left|a_{p_{k} }^{1} \right|e^{-2\pi
i\rho (k)/4} \Delta (\log p_{k} )\right]>>e^{-(1-\delta _{0} )x_{j} } e^{(1-\lambda )x_{j} } >>e^{\delta _{0} x_{j} /2} . \]
Similarly, we can prove an inequality
\[-\sum _{\begin{array}{l} {\quad \quad p_{k} \in P_{3} } \\ {\, \log p_{k} \in \gamma _{j} } \end{array}}Re\left[\left|a_{p_{k} }^{1} \right|e^{-2\pi
i\rho (k)/4} \Delta (\log p_{k} )\right]>>e^{-(1-\delta _{0} )x_{j} } e^{(1-\lambda )x_{j} } >>e^{\delta _{0} x_{j} /2}  .\]
Now fix sequence of values  of $h$: $h_{1},h_{2},...$  taking $h_{1}$ is sufficiently large and writing for    $r>1$ $h_{r}$=$h_(r-1) (1+loq^{-1}h_{r-1})$. It is clear that for each interval $\tau_j$ contains an interval $[h,h(1+log^{-1}h)]\subset[e^{\alpha},e^{\beta}]$ at some $h=h_{r}$.

As it was noted above, the condition (8) is satisfied for unbounded sequence $\left(u_{j}\right)$. Therefore, the intervals  $\tau_{j}$  can be taken disjoint. Thus, from the found above one deduces the existence of an infinite set of indexes $j$ satisfying the conditions
\[\sum _{\log p_{k} \in \tau _{j} ,p_{k} \in P_{v} }(\eta _{k} (s),\varphi (s)) >>e^{\delta _{0} x_{j} /2} ;\, v=0,2,\]
and an infinite set of other values $j$, for which
\[-\sum _{\log p_{k} \in \tau _{j} ,p_{k} \in P_{v} }(\eta _{k} (s),\varphi (s)) >>e^{\delta _{0} x_{j} /2} ;\, v=1,3.\]
Further, from proved above estimations we conclude that
\[|\Delta (x)|\le \pi R^{2} e^{-x/2} .\]
So, $|(\eta _{k} (s),\varphi (s))|\to 0$ when $k\to \infty $. Hence, the series
\[\sum _{n=1}^{\infty } (\eta _{k} (s),\varphi (s))\]
contains two subseries diverging, accordingly, to $+\infty $ and to $-\infty $ having not the common components. Then, some permutation of the series
\[\sum _{n=1}^{\infty } (\eta _{k} (s),\varphi (s))\]
converges conditionally. Therefore, by the theorem 1, \S 6 of [4] there is a permutation of the series $\sum _{p_{n} >y} u_{n} (s)$ converging to $\varphi (s)-\sum _{p_{n} \le y} u_{n} (s)$ uniformly. Taking a long enough partial sum, we receive the necessary result. The lemma 3 is proved.

\textbf{\textit{Note 1.}} The statement of the lemma 3 remains invariable if we consider instead of the disc $\left|s-3/4\right|\le r<1/4$ any disc of a kind $\left|s-\sigma _{0} \right|\le r<r_{0} ;$ $1/2<\sigma _{0} <1$.

\section{The basic auxiliary result.}

Let $\omega \in \Omega $, $\Sigma (\omega )=\{ \sigma \omega |\sigma \in \Sigma \} $ and $\Sigma '(\omega )$ designates the closed set of all limit points of the sequence $\Sigma (\omega )$. For real $t$ we denote $\{ t\Lambda \} $$=(\{ t\lambda _{n} \} )$, where $\Lambda =(\lambda _{n} )$. We will assume that $\mu $ designates the product of the linear Lebesgue measures in $[0,1]$:$\mu =m\times m\times \cdot \cdot \cdot $. In the set $\Omega $ it is possible to define Tychonoff's metric by following equality
\[d(x,y)=\sum _{n=1}^{\infty }e^{1-n} \left|x_{n} -y_{n} \right| . \]

\textbf{\textit{Note}}\textit{ }\textbf{\textit{2}}\textit{.} For completeness we shall show that in the cube $\Omega $ a regular measure may be constructed by using of open sets. At first, we define the volume of the disc of a radius $r>0$:
\[B(0,r)=\left\{E\in \Omega |d(x,0)<r\right\}.\]
Since $|x_{n} |\le 1$ then for the natural number $N$ we have
\[\sum _{n=N+1}^{\infty }e^{1-n} |x_{n} |\le e^{-N} \sum _{n=0}^{\infty }e^{-n} <e^{1-N}   .\]
Taking arbitrarily small real number $\varepsilon >0$ we get
\[\sum _{n=1}^{N}e^{1-n} |x_{n} |\le d(x,0)\le \sum _{n=1}^{N}e^{1-n} |x_{n} |+\varepsilon   \]
when $N\ge \log e\varepsilon ^{-1} $. Therefore,when $B_{N} (0,r)$ denotes  the projection of the sphere $B(0,r)$ into the subspace of first \textit{N} coordinate axis then
\[B_{N} (0,r-\varepsilon )\times [0,1]\times \cdots \subset B(0,r)\subset B_{N} (0,r)\times [0,1]\times \cdots \]. Then, for the volume of $\mu_{N}(r)$ of the ball $B_{N} (0,r)$, we have (see [18, p.319])
\[\mu _{N} (r)-\mu _{N} (r-\varepsilon )=\int _{r-\varepsilon \le \sum _{n=1}^{N} e^{1-n} |x_{n} |\le r}\quad dx_{1} \cdots dx_{N}\]
\[=\textit{2N}\int _{r-\varepsilon \le u\le r}\quad du \int _{\sum _{n=1}^{N} e^{1-n} u_{n} =u}\quad \frac{ds}{\left\| \nabla \right\| } \le \]
\[\le \varepsilon 2^{N} \int _{M}\quad \frac{ds}{\left\| \nabla \right\| }  ,\]
and the last integral is an surface integral over the surface \textit{M} defined by the equation
\begin{equation} \label{9}
\sum _{n=1}^{N}e^{1-n} u_{n}  =u,\quad 0\le u_{k} \le 1;
\end{equation}
here $\nabla $ is a gradient of the linear function on the left side of the latest equality, i.e.
\[\left\| \nabla \right\| =\sqrt{1+e^{-2} +\cdots +e^{2-2N} } .\]
Defining \textit{u1} from (9) we get
\[\int _{M}\quad \frac{ds}{\left\| \nabla \right\| } \le \int _{0}^{1}\cdots \int _{0}^{1}du_{2} \cdots du_{N} =1   .\]
So, we have
\begin{equation} \label{10}
\mu _{N} (r)-\mu _{N} (r-\varepsilon )\le \varepsilon 2^{N} .
\end{equation}
By taking the greatest $N$ satisfying the condition $N\ge \log e\varepsilon ^{-1} $, i.å. $N=[\log e\varepsilon ^{-1} ]+1$ we may write $\varepsilon^{} \le e^{2-N} $. Then from (10) it follows that
\[\mu _{N} (r)-\mu _{N} (r-\varepsilon )\le 2^{N} e^{2-N} \to 0\]
as $N\to \infty $, or as $\varepsilon \to 0$. Since the sequence $\left(\mu _{N} (r)\right)$ is monotonically decreasing then
\[B_{N+1} (0,r)\subset B_{N} (0,r)\times [0,1].\]
So, it is bounded with the lower bound $\mu _{N_{0} } (r/2)$ with $N_{0} $$=[\log 2er^{-1} ]+1$. Therefore, there exists a limit
\[\mathop{\lim }\limits_{\varepsilon \to 0} B_{N} (0,r-\varepsilon )=\mathop{\lim }\limits_{N\to \infty } B_{N} (0,r)=\mu (r)\]
which we receive as a measure of the ball $B(0,r)$.

On this bases it may be introduced the measure in the $\Omega $ by known way by using of open balls. An open ball in the $\Omega $ we define as an intersection $\Omega \bigcap B(\theta ,r)$. The elementary set we define as a set being gotten by finite number of operations of unionize, taking differences or complements of open balls. The outer and inner measures could be introduced by known way (see [17,26]). This measure will be, as it is seen from the reasoning above, a regular measure. As it is clear (see [17, p. 182]), every measurable set in the meaning of introduced measure is measurable in the meaning of product Lebesgue measure also. Main difference between this measure and Haar or product measures is studied in [31]. Now for us it is enough that every set of zero measure can be overlapped by enumerable union of balls with an arbitrarily small total measure.

\textbf{Definition 2.} \textit{Let }$\sigma :N\to N$\textit{ be any one-one mapping of the set of natural numbers. If there will be found a natural number}$m$\textit{ such that }$\sigma (n)=n$\textit{ for any }$n>m$\textit{ then we will say that }$\sigma $\textit{ is a finite permutation. Subset
}$A\subset \Omega $\textit{ we will call  finite-symmetrical if for any element }$\theta =(\theta _{n} )\in A$\textit{ and any finite permutation}$\sigma $\textit{ one has }$\sigma \theta =(\theta _{\sigma (n)} )\in A$\textit{.}

The following lemma is a result of the work [28, p. 46].

\textbf{Lemma 4}. \textit{Let }$A\subset \Omega $\textit{be a finite-symmetrical subset of zero measure and }$\Lambda =(\lambda _{n} )$\textit{ is an unbounded, monotonically increasing sequence of positive real numbers any finite subfamily of elements of which are linearly independent over the field of rational numbers. Let }$B\supset A$\textit{be any open subset in Tychonoff's metric  with }$\mu (B)<\varepsilon $\textit{ and }
\[E_{0} =\{ 0\le t\le 1|\{ t\Lambda \} \in A\wedge \Sigma '\{ t\Lambda \} \subset B\} .\]
\textit{Then, }$m(E_{0} )\le 6c\varepsilon $\textit{ where }$c$\textit{ is an absolute constant and }$m$\textit{ designates the Lebesgue measure.}

\textit{ Proof}. Let $\varepsilon $ be any small positive number. As the numbers $\lambda _{n} $ are linearly independent, for any finite permutation $\sigma $ one has $(\{ t_{1} \lambda _{n} \} )\ne (\{ t_{2} \lambda _{\sigma (n)} \} )$ when $t_{1} \ne t_{2} $. Really, otherwise we would receive the equality $\{ t_{1} \lambda _{s} \} =\{ t_{2} \lambda _{s} \} $ for a great enough natural $s$, i.e. $(t_{1} -t_{2} )\lambda _{s} =k,$ $k\in Z$. Further, writing down the same equality for some other whole $r>s$, we at some whole $k_{1} $ get the relation
\[k_{1} /\lambda _{r} -k/\lambda _{s} =\frac{k_{1} \lambda _{s} -k\lambda _{r} }{\lambda _{r} \lambda _{s} } =0\]
which contradicts the linear independence of the numbers $\lambda _{n} $. Hence, for any pair of various numbers $t_{1} $ and $t_{2} $ one has $(\{ t_{1}
\lambda _{n} \} )\notin \{ (\{ t_{2} \lambda _{\sigma (n)} \} )|\sigma \in \Sigma \} $. By the conditions, there exist a family of open balls $B_{1},B_{2} ,...$ (in Tychonoff's metric) such that each ball does not contain any other ball from this family (the ball containing in other one can be omitted) and
\[A\subset B\subset \bigcup _{j=1}^{\infty } B_{j} ,\sum  \mu (B_{j} )<1.5\varepsilon .\]
Now we take some permutation $\sigma \in \Sigma $ defined by the equalities $\sigma(1)=n_{1} ,$ $...,\sigma (k)=n_{k} $ where natural numbers are taken as below. At first we take $N$ such that
\[\mu (B'_{N} )<2\varepsilon _{1} \]
where $B'_{N} $ is a projection of the ball $B_{1} $ into the subspace of first $N$ co-ordinate axes and $\mu (B_{1} )=\varepsilon _{1} $.  Let $B'_{N} $ be enclosed into the union of cubes with an edge $\delta $ and a total measure not exceeding $3\varepsilon _{1} $. We will put $k=N$ and define
numbers $n_{1} ,...,n_{k} $ using following inequalities
\begin{equation} \label{11}
\lambda _{n_{1} } >1,\lambda _{n_{2} }^{-1} <(1/4)\delta \lambda _{n_{1} }^{-1} ,\lambda _{n_{3} }^{-1} <(1/4)\delta \lambda _{n_{2} }^{-1}
,...,\lambda _{n_{k} }^{-1} <(1/4)\delta \lambda _{n_{k-1} }^{-1} ,\delta <0.1.
\end{equation}
Now we take any cube with an edge $\delta $ and with the center in some point $(\alpha _{m} )_{1\le m\le k} $. Then the point $(\{ t\lambda _{n_{m} } \})$ belongs to this cube if
\begin{equation} \label{12}
|\{ t\lambda _{n_{m} } \} -\alpha _{m} |\le \frac{\delta}{2}.
\end{equation}
Since the interval $(\alpha _{m} -\delta /2,\alpha _{m} +\delta /2)$ has a length $<0.1$ then the real numbers $t\lambda _{n_{m} } $ fractional parts of which lie in this interval have one and the same integral parts during continuous variation of $t$. So at $m=1$ for some whole $r$ one has:
\begin{equation} \label{13}
\frac{r+\alpha _{1} -\delta /2}{\lambda _{n_{1} } } \le t\le \frac{r+\alpha _{1} +\delta /2}{\lambda _{n_{1} } } .
\end{equation}
The measure of a connected set of such $t$ does not exceed the size $\delta \lambda _{n_{1} }^{-1} $. The number of such intervals corresponding to different values of $r=[t\lambda _{n_{1} } ]\le \lambda _{n_{1} } $ does not exceed
\[[\lambda _{n_{1} } ]+2\le \lambda _{n_{1} } +2.\]
So, the total measure of intervals satisfying (12) at $m=1$ is less or equal to
\[(\lambda _{n_{1} } +2)\delta \lambda _{n_{1} }^{-1} \le (1+2\lambda _{n_{1} }^{-1} )\delta .\]
Consider the case $m=2$. Taking one of intervals of a view (12) we will have
\begin{equation} \label{14}
\frac{s+\alpha _{2} -\delta /2}{\lambda _{n_{2} } } \le t\le \frac{s+\alpha _{2} +\delta /2}{\lambda _{n_{2} } } ,
\end{equation}
with some $s=[t\lambda _{n_{2} } ]\le \lambda _{n_{2} } $. As we consider the condition (12) for values $m=1$ and $m=2$ simultaneously, we should estimate a total measure of intervals (14) which have nonempty intersections with intervals of a kind (13) using conditions (11). Every interval of a kind (14) is placed only in one interval with the length $\lambda _{n_{2} }^{-1} $ (on the end points of this interval $t\lambda _{n_{2} } $ takes consecutive integral values) corresponding one and the same value of $s$. The number of intervals with the length $\lambda _{n_{2} }^{-1} $ having a nonempty intersection with one fixed interval of a kind (13) does not exceed the size
\[[\delta \lambda _{n_{1} }^{-1} \lambda _{n_{2} } ]+2\le \delta \lambda _{n_{1} }^{-1} \lambda _{n_{2} } +2.\]
So, the measure of values $t$ for which intervals (14) have a nonempty intersections only with one of intervals of a kind (13) is bounded by the value
$(2+\delta \lambda _{n_{1} }^{-1} \lambda _{n_{2} } )\delta \lambda _{n_{2} }^{-1} $. Since the number of intervals (13) is no more than $\lambda_{n_{1} } +2$, then the measure of a set of values $t$ for which the condition (12) for both numbers $m=1$ and $m=2$ are satisfied simultaneously will be less or equal than
\[(\lambda _{n_{1} } +2)(2+\delta \lambda _{n_{1} }^{-1} \lambda _{n_{2} } )\delta \lambda _{n_{2} }^{-1} .\]
It is possible to continue these reasoning considering all of conditions of a kind
\[\frac{l+\alpha -\delta /2}{\lambda _{n_{m} } } \le t\le \frac{l+\alpha +\delta /2}{\lambda _{n_{m} } } ,m=1,...,k.\]
Then we find the following estimation for the measure $m(\delta )$ of a set of those $t$ for which the points $(\{ t\lambda _{n_{m} } \} )$ located in the given cube with the edge $\delta $:
\[m(\delta )\le (2+\lambda _{n_{1} } )(2+\delta \lambda _{n_{1} }^{-1} \lambda _{n_{2} } )\cdots (2+\delta \lambda _{n_{k-1} }^{-1} \lambda _{n_{k} }
)\delta \lambda _{n_{k} }^{-1} \le \delta ^{k} \prod _{m=1}^{\infty }(1+2m^{-2} ). \]
Summarizing over all such cubes we receive the final estimation of a kind $\le 3c\varepsilon _{1} $ for the  measure of a set of those $t$ for which $(\{ t\lambda _{n_{m} } \} )\in B_{1} $ with an absolute constant $c=\prod _{m=1}^{\infty }(1+2m^{-2} ) $.

We notice that the sequence $\Lambda =(\lambda _{n} )$ satisfying the conditions (11) defined above depends
on $\delta $. We, for each ball $B_{k} $, will fix some sequence$\Lambda _{k} $ using conditions (11). Considering all such balls we designate $\Delta _{0} =\{ \Lambda _{k} |k=1,2,...\} $.

Let's prove that for any point $t\in E_{0} $ the set $\Sigma (\{ t\Lambda \} )$ is contained in the finite union $\bigcup _{k\le n}B_{k} $ for some $n$. Really, let at some $t\in E_{0} $ all members of the sequence $\Sigma (\{ t\Lambda \} )$ does not contained in the union
$\bigcup _{k\le n}B_{k}  $ for any natural $n$. Two cases are possible: 1) there will be found a point $\bar{\theta }\in \Sigma (\{ t\Lambda \} )$ belonging to infinite number of balls $B_{k} $; 2) there will be found a sequence of elements $\bar{\theta }_{j} ,\bar{\theta }_{j} \in \Sigma (\{ t\Lambda \} )$
which does not contained in any finite union of balls $B_{k} $. We will consider both possibilities separately and will prove that they lead to the contradiction.

1) Let $\bar{\theta }\in B_{k_{1} } $,$B_{k_{2} } ,B_{k_{3} } ,...$ are all balls to which the element $\bar{\theta }$ belongs. We shall denote $d$ the distance from $\bar{\theta }$ to the bound of $B_{k_{1} } $. As $B_{k_{1} } $ is an open set, then $d>0$. Let $B_{k} $ be any ball of radius $<d/2$ from the list above containing the point $\bar{\theta }$. From the told it follows that the ball $B_{k} $ should contained in the ball $B_{k_{1}}$. But it contradicts the agreement accepted above.

2) Let $\bar{\theta }$ be some limit point of the sequence ($\bar{\theta }_{j} $). According to the condition of the lemma 3 $\bar{\theta }\in B_{s} $ for some $s$. Let $d$ denote the distance from $\bar{\theta }$ to the bound of $B_{s} $. As $\bar{\theta }$ is a limit point then a ball with the
center in the point $\bar{\theta }$ and radius $d/4$ contains an infinite set of members of the sequence ($\bar{\theta }_{j} $), say members $\bar{\theta }_{j_{1} } ,\bar{\theta }_{j_{2} } ,...$. According to 1) each point of this sequence can belong only to finite number of balls. So, the specified sequence will be contained in a union of infinite subfamily of balls $B_{k} $. Among them will be found infinitely many number of balls having radius $<d/4$. All of them, then, should contained in the ball $B_{s} $. The received contradiction excludes the case 2) also.

So, for any $t\in E_{0} $ it will be found such $n$ for which $\Sigma (\{ t\Lambda \} )\subset $$\bigcup _{k\le n}B_{k}  $. From here it follows that
the set $E_{0} $ can be represented as a union of subsets $E_{k} ,k=1,2,{\rm ...}$, where
\[E_{k} =\{ t\in E_{0} |\Sigma (t\Lambda )\subset \bigcup _{s\le k} B_{s} \} .\]

Therefore,
\[E_{0} =\bigcup _{k=1}^{\infty } E_{k} ;\; \; E_{k} \subset E_{k+1} (k\ge 1).\]
Further, $m(E_{0} )=\mathop{\lim }\limits_{k\to \infty } m(E_{k} )$ in agree with [42, p. 368]. As the set $E_{k} $ is a finite symmetrical, then the
measure of a set of values $t$, interesting us, is possible to estimate by using of any sequence $\Lambda _{k} $, since, as it has been shown above,
the sets $\Sigma (\{ t\Lambda \} )$ for different values of $t$ have empty intersection. So,
\[m(E_{k} )\le \mathop{\lim \sup }\limits_{\Lambda '\in \Delta _{0} } m(E_{k} (\Lambda ')),\]
where $E_{k} (\Lambda ')=\{ t\in E_{k} |(\{ t\Lambda '\} )\in \bigcup _{s\le k}B_{s}  \} $. Hence,
\[m(E_{k} (\Lambda '))\le \sum _{s\le k}m (E^{(s)} (\Lambda ')),\]
where$E^{(s)} (\Lambda ')=\{ t\in E_{0} |(\{ t\Lambda '\} )\in B_{s} \} $. Applying the inequality found above, we receive:
\[m(E_{k} (\Lambda '))\le 6c(\varepsilon _{1} +\cdots +\varepsilon _{k} ).\]
This result invariable for all $\Lambda '=\Lambda _{r} $ beginning from some natural $r=r(k)$. Taking limsup as $k\to \infty $ we receive the
demanded result. The proof of the lemma 4 is finished.

\section{Local approximation.}

\textbf{Lemma 5}. \textit{Let the conditions of the theorem be executed. Then there exist sequences of  points }$\left(\theta _{k} \right)$
\textit{(}$\theta _{k} \in \Omega $\textit{) and natural numbers $\left(m_{k} \right)$ such that}
\[\mathop{\lim }\limits_{k\to \infty } F_{k} \left(\sigma _{0} +s,\; \theta _{k} \right)=F\left(s_{0} +s\right)\, ,s_{0} =\sigma _{0} +it_{0} \]
\textit{uniformly by }$s$\textit{ in the disc }$\left|s\right|\le r<r_{0} $\textit{. }

\textit{Proof. }Let $y>2$ be a whole positive number which will be precisely defined below. We set
\[y_{0} =y,\; y_{1} =2y_{0} ,...,y_{m} =2y_{m-1} =2^{m} y_{0} ,....\]
From the lemma 2 it follows that for given $\varepsilon $ and an integer $y>2$ there exist a finite set $M_{1} $ of primes  such that $M_{1} $ contains all of  prime numbers $p,\, p\le y$ and
\[\mathop{\max }\limits_{\left|s\right|\le r} \left|F(s_{0} +s)-\eta _{1} (s_{1} )\right|\le \varepsilon ;\; \eta _{1} (s_{1} )=\prod _{p\in M_{1}
}f_{p} (e^{-2\pi i(\theta _{p}^{0} +\gamma _{p} )} p^{-s_{1} } ),\; s_{1} =\sigma _{0} +s; \]
besides, $\theta _{p}^{0} =0$ and $\gamma _{p} =(t_{0} /2\pi )\log p$ when $p\le y$, and $\gamma _{p} =0$ if $p>y$. Now we designate
\[h_{1} (s_{1} ;\theta )=F_{1} (s_{1} ;\theta )\cdot \eta _{1}^{-1} (s_{1} )-1, \]
where
\[F_{1} (s_{1} ;\theta )=\prod _{p\le m_{1} }f_{p} \left(e^{-2\pi i(\theta _{p} +\gamma _{p} )} p^{-s_{1} } \right)\, ;\;  \]
$\theta _{p} =\theta _{p}^{0} $ when $p\in M_{1} $ and $m_{1} =\mathop{\max }\limits_{m\in M_{1} } m$. If $r+\delta +2\lambda <r_{0} $ then
\[\int _{\Omega _{1} }\left(\iint \nolimits _{\left|s\right|\le r+\delta +\lambda }\left|h_{1} (s_{1} ;\theta )\right|^{2} d\sigma dt \right) d\theta
\le \iint \nolimits _{\left|s\right|\le r+\delta +\lambda }\left(\int _{\Omega _{1} }\left|h_{1} (s_{1} ;\theta \right|^{2} d\theta  \right) d\sigma
dt\le \]
\[\le \pi (r+\delta +\lambda )^{2} \mathop{\max }\limits_{\left|s\right|\le r+\delta +\lambda } \int _{\Omega _{1} }\left|\sum _{n>y}a_{n} (\theta )
n^{-s_{1} -it_{0} } \right| ^{2} d\theta ;\]
here the summation under the sign of integral is taken over a set of such natural numbers $n$ the canonical factorization of which contains only
primes $p$ with the conditions $p\notin M_{1} ,p\le m_{1} $:
\[a_{n} (\theta )=\prod _{p|n}\left|a_{p}^{\alpha _{p} } \right| e^{2\pi i\alpha _{p} \theta _{p} } ;\; n=\prod p^{\alpha _{p} }  ,\]
and $\Omega _{1} $ means a projection of $\Omega $ into the subspace of co-ordinate axes $\theta _{p} $, $p\notin M_{1} $. By using of orthogonality of the system of functions$e^{2\pi ir\theta } $, $r=1,2,...$ we get
\[\int _{\Omega _{1} }\left(\iint \nolimits _{\left|s\right|\le r+\delta +\lambda }\left|h_{1} (s_{1} ;\theta )\right|^{2} d\sigma dt \right) d\theta
\le \pi \left(r+\delta +\lambda \right)^{2} \sum _{n>y}\left|a_{n} \right| ^{2} n^{2r+2\delta +2\lambda -2\sigma _{0} } \le \]
\[\le \frac{4c^{2} (\lambda )(r+\delta +\lambda )^{2} }{1-2\sigma _{0} -2r-2\delta } y^{1+4\lambda +2r-2\sigma _{0} +2\delta } .\]
Then, there will be found a point $\theta '_{1} =\left(\theta _{p} \right)_{p\notin M_{1} } $ such that
\[\iint \nolimits _{\left|s\right|\le r+\delta +\lambda }\left|h_{1} (s_{1} ;\theta '_{1} )\right|^{2} d\sigma dt \le \frac{4c^{2} (\lambda )(r+\delta
+\lambda )^{2} }{1-2\sigma _{0} -2r-2\delta } y^{1+4\lambda +2r-2\sigma _{0} +2\delta } ,\]
or
\[\mathop{\max }\limits_{\left|s\right|\le r} \left|h_{1} (s_{1} ;\theta '_{1} )\right|\le \sqrt{2} (\delta +\lambda )^{-1} \left(\frac{1}{2\pi } \iint \nolimits _{\left|s\right|\le r+\delta +\lambda }\left|h_{1} (s_{1} ;\theta '_{1} )\right|^{2} d\sigma dt \right)^{1/2}\]
\[ \le c_{1} (\delta ,\lambda)y^{1/2+2\lambda +\delta +r-\sigma _{0} } ,\]
(see [19, p. 345]) with a constant $c_{1} (\delta ,\lambda )>0$. Then, designating $\theta _{1} =(\theta _{0} ,\theta '_{1} ),\, \theta _{0}
=\left(\theta _{p}^{0} \right)_{p\in M_{1} } $, we will have
\[\mathop{\max }\limits_{\left|s\right|\le r} \left\{\left|F(s_{1} +it_{0} )-F_{1} (s_{1} ;\theta _{1} )\right|\right\}\]
\[\le \mathop{\max
}\limits_{\left|s\right|\le r} \left\{\left|F(s_{1} +it_{0} )-\eta _{1} (s_{1} )\right|+\left|\eta _{1} (s_{1} )\right|\cdot \left|h_{1} (s_{1}
;\theta'_{1} )\right|\right\}\le \]
\[\le \varepsilon +(A+1)c_{1} (\delta ,\lambda )y_{0} ^{1/2+r+2\lambda +\delta -\sigma _{0} } ;\quad y_{0} =y,\]
only if $y_{0} $ satisfies the condition
\[(A+1)c_{1} (\delta ,\lambda )y_{0} ^{1/2+r+2\lambda +\delta -\sigma _{0} } \le \varepsilon ;A=\mathop{\max }\limits_{\left|s\right|\le r}
\left|F(s_{1} +it_{0} )\right|.\]

We replace now $\varepsilon $ by $\varepsilon /2$. There is a finite set of primes $M_{2} $ containing the all of prime  numbers $\le 2y_{0} =y_{1} $
and satisfying according to the lemma 3 an inequality
\[\mathop{\max }\limits_{\left|s\right|\le r} \left|F(s_{1} +it_{0} )-\eta _{2} (s_{1} )\right|\le \varepsilon /2;\]
here
\[\; \eta _{2} (s_{1} )=\prod _{p\in M_{2} }f_{p} (e^{-2\pi i(\theta _{p}^{1} +\gamma _{p} )} p^{-s_{1} } ) , \]
$\theta _{p}^{1} =0$$\wedge $$\gamma _{p} =(t_{0} /2\pi )\log p$ when $p\le y_{1} $, and $\gamma _{p} =0$ when $p>y_{1} $. Similarly to performed
above, we find $\theta '_{2} \in \Omega _{2} $ (here $\Omega _{2} $ is a projection of $\Omega $ into the subspace of coordinate axes $\theta _{p} ,\,
p\notin M_{2} $) such that
\[\mathop{\max }\limits_{\left|s\right|\le r} \left|F(s_{1} +it_{0} )-F_{2} (s_{1} ;\theta _{2} )\right|\le \varepsilon ;\; \theta _{2} =(\theta _{1}
,\theta '_{2} ).\]
Really,
\[\left|F_{2} (s_{1} ;\theta )-\eta _{2} (s_{1} )\right|=\left|\eta _{2} (s_{1} )\right|\cdot \left|h_{2} (s_{1} ;\theta )\right|; h_{1} (s_{2}
;\theta )=F_{2} (s_{1} ;\theta )\cdot \eta _{2}^{-1} (s_{1} )-1.\]
Now taking mean values, we receive
\[\mathop{\max }\limits_{\left|s\right|\le r} \left|h_{2} (s_{1} ;\theta '_{2} )\right|\le \sqrt{2} (\delta +\lambda )^{-1} \left(\frac{1}{2\pi }
\iint
\nolimits _{\left|s\right|\le r+\delta +\lambda }\left|h_{2} (s_{1} ;\theta '_{2} )\right|^{2} d\sigma dt \right)^{1/2} \]
\[\le c_{1} (\delta ,\lambda
)(2y)_{0} ^{1/2+r+2\lambda +\delta -\sigma _{0} } .\]
Therefore,
\[\mathop{\max }\limits_{\left|s\right|\le r} \left\{\left|F(s_{1} +it_{0} )-F_{2} (s_{1} ;\theta _{1} )\right|\right\}\le \]
\[\mathop{\max
}\limits_{\left|s\right|\le r} \left\{\left|F(s_{1} +it_{0} )-\eta _{2} (s_{1} )\right|+\left|\eta _{2} (s_{1} )\right|\cdot \left|h_{2} (s_{1}
;\theta'_{2} )\right|\right\}\le \]
\[\le \varepsilon /2+(A+1)c_{1} (\delta ,\lambda )(2y_{o} )^{1/2+r+2\lambda +\delta -\sigma _{0} } \le 2\cdot 2^{1/2+r+2\lambda +\delta -\sigma _{0} }
\varepsilon ;\; \theta _{2} =(\theta _{1} ,\theta '_{2} ).\]
Repeating similar reasoning for every $k>1$ one finds $\theta _{k+1} =\left(\theta _{k} ,\theta '_{k+1} \right)\in \Omega ,\, $with $\theta _{k}
=\left(\theta _{p}^{k} \right)_{p\in M_{k+1} } $ such that $\theta _{p}^{k} =0$$\wedge $$\gamma _{p} =(t_{0} /2\pi )\log p$ when $p\le y_{k} $, and
$\gamma _{p} =0$ at $p>y_{k} $ for which
\[\mathop{\max }\limits_{\left|s\right|\le r} \left|F(s_{1} +it_{0} )-F_{k+1} (s_{1} ;\theta _{k+1} )\right|\le 2^{1+k(1/2+r+2\lambda +\delta -\sigma
_{0} )} \varepsilon ,\]
where
\[F_{k+1} (s_{1} ;\theta )=\prod _{p\le m_{k+1} }f_{p} \left(e^{-2\pi i(\theta _{p} +\gamma _{p} )} p^{-s_{1} } \right)\, ;\;  \; m_{k+1}
=\mathop{\max
}\limits_{m\in M_{k+1} } m.\]
Since $1/2+r+2\lambda +\delta -\sigma _{0} <0$ then uniformly by $s$, $\left|s\right|\le r$
\[\mathop{\lim }\limits_{k\to \infty } F_{k} (s_{1} ;\theta _{k} )=F(s_{1} +it_{0} ).\]
The lemma 5 is proved.

\section{Proof of the theorem.}

 On the theorem's conditions there exists a real $t_{0} $ such that the function $F(s+it_{0} )$ has not zeros in the disc $\left|s-\sigma _{0}
 \right|\le r<r_{0} =\min (1-\sigma _{0} ,\sigma _{0} -1/2)$ at some $1/2<\sigma _{0} <1$ (in the notations of the lemma 5 $r+\delta +2\lambda <r_{0}
 $). Now we will consider the integrals
\[B_{k} =\int _{\Omega }\left(\iint \nolimits _{\left|s\right|\le r}\left|F_{k+1} \left(s_{1} ;\theta _{k+1} +\theta \right)-F_{k} \left(s_{1} ;\;
\theta _{k} +\theta \right)\right| d\sigma dt\right) d\theta ,\]
where $k=0,\, 1,\, \ldots $, where we accept $F_{0} \left(s_{1} ,\; \theta _{0} +\theta \right)=0$, if $k=0$. Applying Schwartz's inequality and
changing the order of the integration, we find as above (denote $\rho =\pi \left(2^{k} y\right)$):
\[B_{k}^{2} \le 4\pi r^{2} \iint \nolimits _{\left|s\right|\le r}d\sigma d\tau \int _{[0,1]^{\rho } }\left|\prod _{p\le 2^{k-1} y_{0} }f_{p} (e^{-2\pi
i\gamma _{p} } p^{-s-it} ) \right|  ^{2} \prod _{p\le 2^{k-1} y_{0} }d\theta _{p}  \times \]
\[\left(\prod _{2^{k-1} y<p\le 2^{k} y}f_{p} \left(e^{-2\pi i(\theta _{p}^{k+1} +\theta _{p} )} \right) -1\right)\prod _{2^{k-1} y<p\le 2^{k}
y}d\theta_{p}  \le \]
\[c\left(\lambda ,\delta \right)\sum _{n>2^{k-1} y_{0} }n^{4\lambda +2r+2\delta -2\sigma _{0} }  \le \]
\[\le c\left(\lambda ,\delta \right)\left(2^{k-1} y_{0} \right)^{1+4\lambda +2r+2\delta -2\sigma _{0} } \, ;\, c\left(\lambda ,\sigma \right)>0..\]
As $1+r+\delta +2\lambda -\sigma _{0} <0$, then from this estimation it follows the convergence of the series below almost everywhere (i.e. for all
$\theta \in \Omega _{0} $, where $\Omega _{0} $ is a subset of  full measure, and set $A=\Omega \backslash \Omega _{0} $ is
finite-symmetrical):\textit{ }
\[\sum _{k=1}^{\infty }\iint \nolimits _{\left|s\right|\le r}\left|F_{k} \left(s+\sigma _{0} ,\; \theta _{k} +\theta \right)-F_{k-1} \left(s+\sigma
_{0} ,\; \theta _{k-1} +\theta \right)\right|d\sigma d\tau ;\; s=\sigma +i\tau   . \]
According to Yegorov's theorem (see [40, p. 166]) this series converges uniformly in the outside of some open set $\Omega \left(\varepsilon
\right),\mu \left(\Omega \left(\varepsilon \right)\right)\le \varepsilon $ for every given $\varepsilon >0$. Put $\Omega '_{1} =\bigcap _{\varepsilon }\Omega \left(\varepsilon \right) $ we can assume that $\mu \left(\Omega '_{1} \right)=0$, and the set $A\bigcup \Omega '_{1} $ is finite-symmetrical
(otherwise it is possible to take the set of all finite permutations of all its elements). There will be found some countable family of balls
$B_{r}$ with a total measure not exceeding $\varepsilon $ the union of which contains the set $A\bigcup \Omega '_{1} $. For every natural \textit{n} we define the set $\Sigma '_{n} (t\Lambda )$ as a set of all limit points of the sequence $\Sigma _{n} (\bar{\omega })=\{ \sigma \bar{\omega }|\sigma \in
\Sigma \wedge \sigma (1)=1\wedge \cdots \wedge \sigma (n)=n\} $. Let
\[B^{(n)} =\{ t|\{ t\Lambda \} \in A\wedge \sum '_{n} (\{ t\Lambda \} )\subset \bigcup _{r=1}^{\infty }B_{r}  \} ,\lambda _{n} =(1/2\pi )\log p_{n}
n=1,2,....\]
For every $t$ the sequence $\sum _{n+1} (\{ t\Lambda \} )$ is a subsequence of the sequence $\sum _{n} (\{ t\Lambda \} )$. Therefore, $\sum '_{n+1}
(\{t\Lambda \} )\subset $$\sum '_{n} (\{ t\Lambda \} )$ and we have $B^{(n)} \subset B^{(n+1)} $. Then we have an inequality  $m(B)\le \mathop{\sup
}\limits_{n} m(B^{(n)} )$ denoting $B=\bigcup _{n} B^{(n)} $.

Let's estimate $m(B^{(n)} )$. The set $\sum '_{n} (\{ t\Lambda \} )$ is a closed set. Clearly, if we will "truncate" sequences $\{ t\Lambda \} $
leaving only components $\{ t\lambda _{n} \} $ with indexes greater than $n$ and will denote the truncated sequence as $\{ t\Lambda \} '\in \Omega $,
then the set $\sum '(\{ t\Lambda \} ')$ also will be closed. Now we consider the products $[0,1]^{n} \times \{ \{ t\Lambda \} '\} $ (external brackets
designate the set of one element) for every $t$. We have
\[\{ t\Lambda \} \in [0,1]^{n} \times \{ \{ t\Lambda \} '\} \subset A.\]
(The example below shows that from the feasibility of the last relationship it does not follow the equality $A=\Omega $. Let
$I=[0,1];U=[0;1/2];V=[1/2;1]$ and
\[X_{0} =U\times U\times \ldots ,X_{1} =V\times U\times \ldots ,\]
\[X_{2} =I\times V\times U\times \ldots ,X_{s+1} =I^{s} \times V\times U\times \ldots ,....\]
Clearly, that $\mu (X_{s} )=0$for all$s$. Let
\[X=\bigcup _{s=0}^{\infty } X_{s} .\]
So, we have $X=[0,1]^{s} \times X$ for any natural$s$. Then $\mu (X)=0$ and $X\ne \Omega $).
Let $(\theta _{1} ,...,\theta _{n} )\in [0,1]^{n} $. There exist a neighborhood $V\subset [0,1]^{n} $ of this point such that  $(\theta _{1}
,...,\theta _{n} ,\{ t\Lambda \} ')\in V\times W\subset \bigcup _{r}B_{r}  ,\; $for some neighborhood $W$ of the point $\{ t\Lambda \} '$. Since the
set  $[0,1]^{n} $ is closed, then they can be found a finite number of open sets $V$ the union of which contain $[0,1]^{n} $. The intersection of
corresponding open sets $W$ being an open set contains the point $\{ t\Lambda \} '$. Therefore, we have
\[[0,1]^{n} \times \{ \{ t\Lambda \} '\} \subset \bigcup V \times \bigcap W =[0,1]^{n} \times \bigcap W \subset \bigcup _{r\in R} B_{r} ,\]
for each considered point $t$. The similar relationship is fair in the case when the point $\{ t\Lambda \} $ would be replaced by any limit point
$\bar{\omega }$ of the sequence $\Sigma (\{ t\Lambda \} )$ also, because $\bar{\omega }\in B_{r} $. If one denotes by $B'$ the union of all open sets
of a kind $\bigcap _{r\in R} B'_{r} $, corresponding to every possible values of $t$ and of a limit point $\bar{\omega }$, then we will receive the
relation
\[\{ t\Lambda \} \in [0,1]^{n} \times \{ \{ t\Lambda \} '\} \subset A\subset [0,1]^{n} \times B'\subset \bigcup _{r=1}^{\infty }B_{r}  ,\]
for each considered values of $t$ and
\[\{ \bar{\omega }\} \in [0,1]^{n} \times \{ \bar{\omega }\} '\subset A\subset [0,1]^{n} \times B'\subset \bigcup _{r=1}^{\infty }B_{r}  ,\]
for each limit point$\bar{\omega }$. From this it follows the inequality $\mu ^{*} (B')\le \varepsilon $, where $\mu ^{*} $ means an outer measure. The set $B'$ is open and $\Sigma '(\{ t\Lambda \} ')\in B'$. Now we can apply the lemma 3 and receive an estimation $m(B^{(n)} )\le 6c\varepsilon $. Thus,
we have $m(B)\le 6c\varepsilon $.

 Let $t\notin B$. Then, $t\notin B^{(n)} $ for every $n=y_{k} ,k=1,2,3,...$. Consequently, for every $k$ there is a such limit point $\bar{\omega
 }_{k} \in \Omega \backslash \bigcup _{r} B_{r} $ of the sequence $\sum _{n} (\{ t\Lambda \} )$ for which the series
\[\sum _{l=1}^{\infty } \int  \int _{|s|\le r} |F_{l} (\sigma _{0} +s;\bar{\theta }_{l} +\bar{\omega }_{k} )-F_{l-1} (\sigma _{0} +s;\bar{\theta
}_{l-1} +\bar{\omega }_{k} )|d\sigma d\tau \]
converges. As the set $\Omega \backslash \bigcup _{r} B_{r} $ is closed, the limit point $\overline{\omega }=(\{ t\Lambda \} )$ of the sequence
$(\bar{\omega }_{k} )$ will belong to the set $\Omega \backslash \bigcup _{r} B_{r} $. Therefore, the series
\begin{equation} \label{15}
\sum _{l=1}^{\infty } \int  \int _{|s|\le r} |F_{l} (s+\sigma _{0} ;\theta _{l} +i\{ t\Lambda \} )-F_{l-1} (s+\sigma _{0} ;\theta _{l-1} +i\{ t\Lambda
\} )|d\sigma d\tau
\end{equation}
converges. So the last series converges for all $t$ with exception of values $t$ from some set of  measure not exceeding $12c\varepsilon $. Owing to
randomness of $\varepsilon $ the last result shows a convergence of (15) for almost all $t$ (clearly, that the condition $0\le t\le 1$ can be omitted now). Then, by the lemma 1 for $\delta _{0} <1$ taken arbitrarily the sequence
\begin{equation} \label{16}
F_{k} (s+\sigma _{0} ;\theta _{k} +i\{ t\Lambda \} )
\end{equation}
converges, for all such $t$, in the disc $|s|\le r\delta _{0} (\delta _{0} <1)$ uniformly to some analytical function $f(s+\sigma _{0} ;t)$:
\[\mathop{\lim }\limits_{k\to \infty } F_{k} (s+\sigma _{0} +it;\theta _{k} )=f(s+\sigma _{0} ;t).\]
Despite the received result, we cannot use $t$ as a variable as the left and right parts of this equality can differ each from other by their
arguments (the right part is defined as a limit of the sequence (16). where $t$ enters into the expression containing discontinuous function). Hence, the principle of analytical continuation cannot be applied. To finish the theorem's proof we take any large real number $T$. As considered values $t$ are everywhere dense in the segment $[-T,T]$, the union of discs $C(t)=\{ \sigma _{0} +it+s:|s|\le r\delta _{0} \} $ contains the rectangle
\[\sigma _{0} -r\delta _{0}^{2} \le Re(s+\sigma _{0} )\le \sigma _{0} +r\delta _{0}^{2} ,-T\le Im(s+3/4)\le T\]
in which conditions of  the lemma 1 are executed for the series
\begin{equation} \label{17}
F_{1} (s+\sigma _{0} ;\theta _{1} )+(F_{2} (s+\sigma _{0} ;\theta _{2} )-F_{1} (s+\sigma _{0} ;\theta _{1} ))+\ldots .
\end{equation}
Hence, by the lemma 1, this series defines an analytical function in the considered rectangle which coincides with $F(s_{0} +s)$ in the disc
$C(0)$. To apply the principle of analytical continuation we take one-connected open domain where both of the functions $\log F_{*} (s)$ and $\log
F(s+s_{0} )$ are regular (here function $F_{*} (s)$ is the sum of the series (17)). Let $\rho _{1} ,...,\rho _{L} $ designate all possible zeros of
the function $F(s_{0} +s)$ in the considered rectangle the contour of which does not contain zeros of the function $F(s_{0} +s)$. We will take cuts
through the segments $1/2\le Res\le Re\rho _{l} ,$ $Ims=Im\rho _{l} ,$$l=1,...,L.$ In the open domain of the considered rectangle not containing
specified segments the functions $\log F_{*} (s)$ and $\log F(s+s_{0} )$ are regular. Therefore, in this domain the equality $\log F_{*} (s)$$=\log
F(s+s_{0} )$ holds. Then, the equality $F_{*} (s)=F(s+s_{0} )$ is executed in all open domain defined above. Now we receive a justice of the relation
$F_{*} (s)=F(s+s_{0} )$ in the all rectangle (without cuts) where both functions are regular. The theorem is proved.

\section{Proof of the consequence.}

The conclusion of the consequence based on the theorem of Rouch'e ( see [19,  p. 137]). Let $t$ be any real number. We shall prove that for any
$0<r'<3/4$ in the domain bounded by the circle $C'=\{ s||s-\sigma _{0} -it|=r'\} $ the function $F(s)$  has not zeroes. Since there are only a
finite set of zeroes satisfying the condition $|s-\sigma _{0} -it|\le r<3/4$, then we may take $r>r'$ such that the disc $C=\{ s||s-\sigma _{0} -it|=r\} $ does not contain zeroes of $F(s)$. Let
\[m=\mathop{\min }\limits_{s\in C} |F(s)|.\]
Since the $C$ is a compact set, clearly $m>0$. By the theorem there exist $n=n(t)$ such that the following inequality is executed on $C$ or in the
disc bounded by $C$:
\[|F(s)-F_{n} (s;\bar{\theta }_{n} )|\le 0.25m.\]
Then, on contour of $C$ the following inequality is true:
\[|F(s)-F_{n} (s;\bar{\theta }_{n} )|<|F(s)|.\]
Then, from the theorem of Rouch'e it follows that the functions $F(s)$ and $F_{n} (s;\bar{\theta }_{n} )$ have an identical number of zeroes inside
$C$. But, the function $F_{n} (s;\bar{\theta }_{n} )$ has not zeroes there. Hence, $F(s)$ also has not zeroes in the open disc bounded by $C$. As $t$ is taken arbitrarily, from the last we conclude that the strip $-r<Re\, s-3/4<r$ for any $0<r<1/4$ is free from the zeroes of the function $F(s)$. Obviously, for any $1/4>\lambda >0$ there exist a segment $[1/2+\lambda +i\tau ,1-\lambda +i\tau ]$ not containing zeros of $F(s)$. This segment can be covered by finite number of discs not containing zeros of $F(s)$. Applying proved above to each of such discs we receive the strips free from the zeros of $F(s)$ the union of which contains the strip $1/2+\lambda <Res<1-\lambda $. As $\lambda $ is any positive number then the statement of the theorem is proved.

\noindent

\noindent

\bibliographystyle{amsplain}

\bibliography{}
\textbf{References}
\noindent

\noindent 1. L.Euler. Introduction to the analyses of infinitesimals. -- Ì. :ONTI, 1936. (rus).

\noindent 2. B. Riemann. On the number of prime numbers not exceeding a given quality.
//Compositions.--Ì. : ÎGIZ, 1948 -- p. 216 -- 224(rus).

\noindent 3. E. C. Titchmarsh. Theory of Riemann Zeta -- function. Ì. : IL, 1953(rus).

\noindent 4. S. M. Voronin, A.A. Karatsuba. The Riemann Zeta--function. M: fiz. mat. lit. ,1994,
 376 p. (rus).

\noindent 5. À.À. Êàðàöóáà. Îñíîâû àíàëèòè÷åñêîé òåîðèè ÷èñåë. Ì. Íàóêà, 1983, (rus).

\noindent 6. H.L. Montgomery. Topics in multiplicative number theory.Ì.,1974(rus).

\noindent 7. H.Davenport. Multiplicative number theory. Ì.$:$Nauka 1971(rus).

\noindent 8. K. Chandrasekharan. Arithmetical functions. Ì. Nauka, 1975., 270 pp. (rus).

\noindent 9. S.M. Voronin. On The distribution of non -- zero values of the Riemann Zeta
function. Labors. MIAS -- 1972 -- v. 128, p. 153-175, (rus).

\noindent 10. S.M. Voronin. On an differential independence of \textit{$\zeta$} -- function.
 Reports of AS USSR -- 1973. v. 209, ¹ 6, pp.1264 -- 1266,
(rus).

\noindent 11. S.M. Voronin. On an differential independence of Dirichlet`s L-- functions.
 Àctà  Àrith. -- 1975 , v. ÕÕVII -- pp. 493 -- 509.

\noindent 12. S.M. Voronin. The theorem on ``universality'' of  the Riemann zeta -- function.
 Bulletin  Acad. Sci. of  USSR mat.ser. -- 1975 -- v. 39,
¹ 3 -- pp. 475 -- 486, (rus).

\noindent 13. S.M. Voronin. On the zeroes of zeta -- functions of quadratic forms. Labors of
MIAS -- 1976 -- v. 142 -- pp. 135 -- 147 (rus).

\noindent 14. S.M. Voronin. Analytical properties of Dirichlet generating functions of
arithmetical objects: Diss.\dots D -- r of fiz. -- mat. sci.
MIAS USSR -- Ì., 1977 -- 90p, (rus).

\noindent 15. S.M. Voronin. On an zeroes of some Dirichlet series, lying on the critical line.
Bull. Acad. Sci. of  USSR mat.ser. -- 1980 -- v. 44 ¹1
-- pp.63-91 (rus).

\noindent 16. S.M. Voronin. On the distribution of zeroes of some Dirichlet series. Labor.
MIAS -- 1984 -- v. 163 -- pp. 74 -- 77, (rus).

\noindent 17. N.Dunford and J.T.Schwartz. Linear operators. Part I: General theory. Ì. PFL,1962,
 896 ð.

\noindent 18. Êóðàíò R.Differential and integral calculus. Ì: Nauka, 1967.

\noindent 19. E. C. Titchmarsh. Theory of function. Ì.: GITTL, 1951.506 pp.(rus).

\noindent 20. W. Rudin. Principles of mathematical analysis. Ì.: Mir,1976.319 pp. (rus).

\noindent 21. Hewitt E. and Ross K. Abstrakt Harmonic Analysis. v.1, Nauka, 1975.

\noindent 22. B. Bagchi. A joint universality theorem for Dirichlet L --functions. Math. Zeit.,
1982,v. 181, p. 319-335.

23. A. Laurinchikas. Limit Theorems for the Riemann Zeta-Function, Kluwer, Dordrecht, 1996.

\noindent 24. Ëàóðèí÷èêàñ À.Ï. On zeros of linear combinations of numbers Äèðèõëå. Lit.
Mat. collec., 1986, v.26, ¹3, p.468-477.

\noindent 25. A.Zigmund. Trigonometrical series. v. 2., Ì: Mir, 1965.

\noindent 26. V. I. Bogachev. Measure Theory. Springer-Verlag Berlin Heidelberg 2007, v. 1-2.

\noindent 27. Dzhabbarov, I.Sh. Mean values of Dirichlet L-functions on short closed intervals of the critical line and their applications. (Russian. English summary) Izv. Akad. Nauk Az. SSR, Ser. Fiz.-Tekh.Mat. Nauk 1988,  No.1, 3-9 (1988).

\noindent 28. Dzhabbarov I. Sh. On Ergodic Hypothesis. Euler International Mathematical
Institute. Topology, Geometry and Dynamics: Rokhlin Memorial. Short abstracts of an international meeting held on January 11-16, 2010. St. Petersburg, 2010, p. 46-48.

\noindent 29. Dzhabbarov I. Sh. Uniform approximation of Dirichlet series by partial products
of Euler type. International conference ``Approximation
theory''. Abstracts (Saint-Petersburg, 6-8 may, 2010), St. Petersburg, 2010, p. 117-119.

\noindent 30. Jabbarov I. Sh. The Riemann Hypothesis. ArXiv:1006.0381v3, 2010.

\noindent 31. Dzhabbarov I. Sh. On the connection between measure and metric in infinite dimensional
space. International conference on Differential
Equations and Dynamical Systems. Abstracts. Suzdal (Russia) July 2-7, 2010, Moscow, 2010, p. 213-214.
\end{document}